\documentclass{demorgan}
\title[Mathematics inherited]{Mathematics discovered, invented, and inherited}
\author{Alexandre Borovik}
\date{2013}
\journal{Selected Passages from Correspondence with Friends 1 no.~4 (2013), 13--28 \hfill \ }
\usepackage{vmargin}
\usepackage{url}
\usepackage{amssymb}
\usepackage{amsmath}
\usepackage{vmargin}
\usepackage{url}
\usepackage{epsfig}
\usepackage{graphicx}
\newcommand{\beq}{\begin{equation}}
\newcommand{\eeq}{\end{equation}}
\newcommand{\bea}{\begin{eqnarray*}}
\newcommand{\eea}{\end{eqnarray*}}
\newcommand{\bt}{\begin{theorem}}
\newcommand{\bd}{\begin{description}}
\newcommand{\ed}{\end{description}}
\newcommand{\bi}{\begin{itemize}}
\newcommand{\ei}{\end{itemize}}
\newcommand{\bq}{\begin{quote}}
\newcommand{\eq}{\end{quote}}

\extraline{\copyright\ 2013 Alexandre V. Borovik
\hfill www.borovik.net/selecta}
\classno{Primary 20D05, 20F11; secondary 01A60}

\volume{1}
\begin{document}

\setcounter{page}{13}
\maketitle

\begin{abstract} \normalsize
The classical platonist/formalist dilemma in philosophy of mathematics can be expressed in lay terms as a deceptively naive question:

\bq
{\normalsize\emph{Is new mathematics discovered or invented?}}\\[-1.5ex]
\eq

\noindent
Using an example from my own mathematical life, I argue that there is also a third way:

\bq
{\normalsize\emph{new mathematics can also be inherited}}\\[-1.5ex]
\eq
and in the process briefly discuss a remarkable paper by W. Burnside of 1900 \cite{burnside00.269}.
\end{abstract}

\

\begin{flushright}

\emph{Ontogeny recapitulates phylogeny}.\\

Ernst Haeckel

\end{flushright}

\large

\section{Introduction}

Bill Kantor asked me how had it happen that I came to possession of a certain algorithm in computational group theory independently from the remarkable work by John Bray \cite{bray00.241}. I wish to make it clear that I do not claim priority; the algorithm  is John Bray's discovery. But I was lucky to be exposed to a rich mathematical tradition and had happen \emph{to inherit} the principal ideas of the algorithm which could be traced back to 1899 (the year when Burnside's paper \cite{burnside00.269} was submitted for publication) but remained forgotten for the good part of 20th century.

Perhaps it is worth telling this story in detail: it seems to be an interesting case study of the roles of tradition and cultural heritage in  mathematics.

\section{Bray's algorithm.}
\label{sec:centralisers}

To explain the algorithm, I quote \emph{verbatim} this entire section from my paper \cite{borovik02.7}.

\begin{quotation}
It is well known that two involution generate a dihedral subgroup of any group,
finite or not. Indeed, if $u$ and $v$ are involutions then
$$
(uv)^u = u^{-1} \cdot uv \cdot u = vu = (uv)^{-1}
$$
and similarly $(uv)^v = (uv)^{-1}$. Hence $u$ and $v$ invert every
element in the cyclic group $\langle uv\rangle$. If the element
$uv$ is of even order then $u$ and $v$ invert the involution
${\rm i}(uv)\in \langle uv \rangle$ and centralise it. If,
however, the element $uv$ has odd order then, by the
Sylow Theorem, the involutions $u$ and $v$ are conjugate by an
element from $\langle uv \rangle$.

This simple observation, due to Richard Brauer, was the starting
point of his programme of classification of finite simple groups
in terms of centralisers of involutions. Remarkably, in the
context of black box groups it can be developed into an efficient
algorithm for constructing black boxes for the centralisers of
involutions.

Let ${X}$ be an arbitrary black box finite group and assume that
$x^E =1$ for all elements $x\in {{X}}$. Write $E = 2^t \cdot r$
with $r$ odd. Let $x$ be a random element in ${X}$. Notice that

\begin{itemize}
\item if $x$ is of odd order, then $x^r=1$ and
$y = x^{(r+1)/2}$ is a square root of $x$:
$$
y^2 =  x^{(r+1)/2} \cdot x^{(r+1)/2} = x^{r+1} = x;
$$

\item if $x$ is of even order then $x^r$ is a $2$-element
and the consecutive squaring of $x^r$ produces the
involution ${\rm i}(x)$ from the cyclic group $\langle x \rangle$.
\end{itemize}

Furthermore, the elements $y = \sqrt{x}$ and ${\rm i}(x)$ can be
found by $O(\log E)$ multiplications.

Let now $i$ be an involution in ${X}$. Construct a random element
$x$ of ${X}$ and consider $z = ii^x$.

\begin{itemize}
\item If $z$ is of odd order and $y = \sqrt{z}$ then
$$
i^y = y^{-1}i y = i yy = iz = i\cdot ii^x = i^x
$$
and $yx^{-1} \in C_X(i)$. We write\footnote{Adrien Deloro and other model theorists colleagues called $\zeta_1$ \emph{the zeta map}, see a quote from Adrien Deloro on Page~\pageref{deloro-quote}.} $ yx^{-1} = \zeta_1(x)$.

\item if $z$ is of even order then ${\rm i}(z)$
lies in the center of the dihedral group $\langle i, i^x\rangle$ and thus ${\rm i}(z) \in C_X(i)$.
We write ${\rm i}(z) = \zeta_0(x)$.

\end{itemize}
Notice that $\zeta_1(x)$ can be computed without knowing the
order $o(x)$ of $x$. One can test whether an element has odd order
by raising it to the odd part $r$ of $E$, and if $o(x)$ is odd
then $x^{(r+1)/2} = x^{(o(x)+1)/2}$.

\medskip

Thus we have a map $\zeta = \zeta_1 \sqcup \zeta_0$ defined by
\begin{eqnarray*}
\zeta: X & \longrightarrow &  C_X(i)\\
x & \mapsto & \left\{ \begin{array}{ll}
\zeta_1(x) = (ii^x)^{(r+1)/2}\cdot x^{-1} & \hbox{ if } o(ii^x) \hbox{ is odd}\\
\zeta_0(x) = {\rm i}(ii^x)  &  \hbox{ if } o(ii^x) \hbox{ is even.}
\end{array}\right.
\end{eqnarray*}

If $c \in C_X(i)$ then
\begin{eqnarray*}
\zeta_0(xc) & =&  {\rm i}(i\cdot i^{xc}) = {\rm i}(i^c\cdot i^{xc}) =
{\rm i}((i\cdot i^x)^c) = {\rm i}(ii^x)^c\\
& = & \zeta_0(x)^c,\\[2ex]
\zeta_1(cx) & = & (ii^{cx})^{(r+1)/2}\cdot x^{-1}c^{-1} =
(ii^{x})^{(r+1)/2}\cdot x^{-1}c^{-1}\\
& = &\zeta_1(x)\cdot c^{-1}.
\end{eqnarray*}

Hence if the elements $x \in {X}$ are uniformly distributed and
independent in $X$ then
\begin{itemize}
\item  the distribution of elements $\zeta_1(x)$ in $C_X(i)$
is invariant under  right multiplication by elements $c \in C_X(i)$,
that is, if $A \subset  C_X(i)$ and $c \in C_X(i)$ is an arbitrary element then
the probabilities $P(\zeta_1(x) \in A)$ and $P(\zeta_1(x) \in Ac)$ coincide.

\item The distribution of involutions $\zeta_0(x)$ is
invariant under the
action of $C_X(i)$ on itself by conjugation, that is,
$$
P(\zeta_1(x) \in A) =P(\zeta_1(x) \in A^c).
$$
\end{itemize}
\end{quotation}

\section{A few words about how I was learning group theory}

I was 14 years old in Summer 1971 when I started to read a Russian translation of Robin Hartshorne's book \emph{Foundations of Projective Geometry} \cite{hartshorne67}; I bought the book for 45 kopecks (the cost of three standards portions of ice cream) in the village book shop. Chapter III of the book contained some basic group theory, and I was mesmerised by its level of abstraction. A few months later, when I moved from my home village in Eastern Siberia to a boarding school in Novosibirsk\footnote{My brief notes on life and study in this school can be found in \cite{borovik12.23}.}, I continued my attempts to learn some group theory using some random Russian textbooks of algebra. My mathematics tutorial class teacher, Sergei Syskin, discovered that and was appalled. ``Why do you read all that rubbish?'' -- he said -- ``if you want to learn group theory, do it properly and read proper books''. Next day he brought me \emph{Abstract Group Theory} by Otto  Schmidt\footnote{In Russia, Otto Yulievich Schmidt was famous mostly for his daring Arctic expedition on \emph{SS Chelyuskin}; he was also a mathematician, a geophysicist, and a high-ranking member of the ruling Communist establishment who at one point was in charge of the entire Russian sector of Arctic. With melting ice opening up Arctic for commercial exploitation, his hame is likely to resurface in the annals of history. The most recent book on him \cite{matveeva93} is still a rather one-sided official biography of the much more interesting and controversial man.} printed in 1933\footnote{An English translation from this edition is \cite{schmidt66}.}, with brittle yellowish pages,  archaic notation, Gothic fonts, but also with an obvious feel of something which was ground-breaking for its time---especially after taking into consideration that the book was republished, with minor changes, from its first edition of 1916\footnote{Schmidt finishes his \emph{Introduction to the Second Edition} with the words:\\[-2ex] \bq{\footnotesize \dots an abundance of other work does not permit me in the immediate future to undertake [\dots] a revision. For the present, therefore, I have accepted the suggestion of the State Technical and Theoretical Publishing House to republish the book in its original form. I have made only insignificant changes.\\[2ex]
August, 1933 \hfill O. Schmidt\\
The Icebreaker \emph{Cheluskin}}\eq}. I spent my school vacations in January 1972 reading this book---there were nothing else to do in my snowbound home village.

At that time Syskin was a young graduate student who however had already published some good papers in finite group theory. His choice of the book was apparently influenced by the educational Zeitgeist of that peculiar place and time: coming back to the origins and to the fundamentals. Or perhaps he was influenced by his own teacher, Victor Danilovich Mazurov: when a few years later Mazurov started a study group (for just three undergraduate students: Elena Bryukhanova, Evgenii Khukhro, and myself), it was focused on celebrated ground-breaking papers, such as Hall and Higman \cite{hall-higman56.1}.\footnote{We met three times a week, each session lasting 90 minutes. I gather Mazurov was spending on us between 6 to 9 hours a week, every week, for several years.}

I think it was at that time in Novosibirsk that I picked up  Ernst Haeckel's maxim
\bq
{\normalsize
``\emph{ontogeny recapitulates phylogeny},''
}
\eq
it was very much part of Zeitgeist.

But let us return to 1972. A few months later Syskin gave me Gorenstein's \emph{Finite Groups} \cite{gorenstein1968}. This book offered a completely different perspective on (finite) groups. Schmidt treated groups elementwise, subgroups being not just subsets but some kind of containers, which had to be filled with elements, otherwise they had no much value. Gorenstein was using a functorial approach; for him, really interesting subgroups were values of \emph{functors}, and in many instances he did not care much about the elements inside. For me, it was a dramatic and revolutionary change---and for that reason, Chapter 9, \emph{Groups of even order}, of Gorenstein's book caught my eye---the exposition of the Brauer-Fowler theory was very old-fashioned and very elementwise. These results and proofs would not be out of place in Schmidt's \emph{Abstract Group Theory}.

\section{William Burnside, 1900}

Much later, in about 1982--83, I attempted to apply methods of the classification of finite simple groups to model-theoretic algebra, more specifically to the study of simple $\omega$-stable groups of finite Morley rank. The origins of this project are explained in the following quote from \emph{Introduction} to  \cite{borovik_ABC_2008}, the book that summarised 25 years of work of my colleagues and myself on this programme:

\bq
{\normalsize
Vladimir Nikanorovich Remeslennikov in 1982 drew my attention
to Gregory Cherlin's paper \cite{cherlin1979.1} on groups of finite Morley
rank and conjectured that some ideas from my work [on periodic
linear groups] could be used in this then new area of algebra.
}
\eq

An attempt to work in the context of Cherlin's paper lead me to a surprising  discovery. It had become apparent that when the statement of the Brauer-Suzuki-Wall Theorem \cite{brauer58.718}, one of the first and basic characterisation of finite simple groups in terms of centralisers of involutions, had been adequately transferred  to the infinite domain, it remained valid not only for groups like ${\rm PSL}_2(\mathbb{C})$, say, but also for ${\rm SO}_3(\mathbb{R})$. The former could be seen as a natural analog of groups ${\rm PSL}_2(q)$ characterised by Brauer, Suzuki, and Wall, but the latter was an anisotropic algebraic group and had no analogues among finite groups. It became clear to me that I had to revisit the early chapters of the classification of finite simple groups and rework them starting from the first principles in the new model theoretic context.

This is why I red the mathematical biography of Richard Brauer written by Walter Feit \cite{feit79.1} with great interest. I was struck by the following paragraph (pp. 11--12).

\bq
{\normalsize
A few years ago, as I was preparing a lecture on the history of group
theory, I came across a paper of Burnside in which he characterized the
groups ${\rm SL}_2(2^a)$ for $a > 1$ as the only simple groups of even order in which the
order of every element is either $2$ or odd \cite{burnside00.269}. In this paper Burnside used
some of the basic properties of involutions in a way quite similar to the way
that Brauer used them fifty years later. However Burnside did not realize the
importance of this approach. His paper had an innocuous title and appeared
in a journal that is not readily available in mathematical libraries. I don't
believe that he referred to the paper in his book \cite{burnside1911} and I don't know of
any mathematical paper by Burnside or anyone else that refers to this paper.
When I told Brauer about this paper, he was as surprised to hear of it as I
had been when I first found it. It is a tribute to Brauer's insight that he
realized that these extremely elementary arguments concerning involutions
are of fundamental importance. The fact that a mathematician of the caliber
of Burnside had overlooked the importance of such arguments, even after
proving the result in \cite{burnside00.269}, shows that this insight was far from obvious.
}
\eq

Burnside's paper was absolutely unavailable to me (after all, I lived in Omsk at that time), and I jumped at the opportunity to recover from first principles his proof of his result. This was easy. I was much helped by three crucial pieces of knowledge:
\bi
\item I knew and loved the Brauer-Fowler theory;
\item comments by Feit provided some powerful hints;
\item I had never read Burnside, but Schmidt's book provided a good insight into the style of group-theoretic thinking of that time---and I needed to learn this style of thinking.
\ei
My proof of Burnside's Theorem used two simple ideas: if $G$ is a finite simple group where centraliser of every involution is an elementary abelian $2$-group, then
\bi
\item[(a)] conjugation of involutions is ``local,'' that is, if $i$ and $j$ are two non-commuting involutions, then they are conjugated in the dihedral group $\langle i,j\rangle$ (and, moreover, are conjugated by some element $y$ from the cyclic group $\langle ij \rangle$ of odd order);
\item[(b)] the normaliser $B = N_G(T)$ of a $2$-Sylow subgroup $T$ of $G$ had to be filled with elements.
\ei
At that time, I was already thinking about transferring the technique to groups of finite Morley rank, and I had to explicitly express the conjugating element $y$ in (a) in terms of $i$ and $j$; to my great relief, I discovered that $y$ could be taken being equal
\[
y = \sqrt{ij}
\]
because $i$ inverts $z = \sqrt{ij}$ and (I repeat the calculation from Section~\ref{sec:centralisers})
\[
y^{-1}iy = iyy = i\sqrt{ij}\sqrt{ij} = i\cdot ij = j,
\]
or, in brief,
\[
i^{\sqrt{ij}} = j.
\]
This calculation was very important for me because it worked in groups of finite Morley rank where the cyclic group $\langle ij\rangle$ could be  embedded in a $2$-divisible abelian group normalised and inverted by $i$. This version of Brauer's famous lemma helped me to prove the conjugacy of $2$-Sylow subgroups in groups of finite Morley rank \cite{borovik90.478}.

Regarding (b), I was using the ``double conjugation'' trick: take involutions $t, s \in T$ and some involution $r \in G\smallsetminus T$ so that $r$ commutes neither with $t$ nor $s$, then if you take
\[
 b= \sqrt{tr}\sqrt{rs}
\]
you have
\[
t^b = (t^{\sqrt{tr}})^{\sqrt{rs}} = s
\]
and since $T = C_G(t) = C_G(s)$, $b$ normalises $T$ hence belongs to  $B$.

I had not even thought about publication of my proof---indeed what was the point to republish a theorem from 1900? But the
tricks (a) and (b),  ``conjugation by a square root'' and ``double conjugation'',  which I learnt from this exercise with help from  Burnside, Brauer, and Feit, allowed me to start exploration of groups of finite Morley rank, see more on that in Section~\ref{sec:Euler}.

\section{But how did Burnside do it?}

As I said before, I was trying to recover Burnside's original proof, but was not able to compare my proof with the original until this year, 2013, when I had finally discovered Burnside's paper on the Internet. Meanwhile a few years ago I encountered  a review MR0352256 (written by Ulrich Dempwolff) of David Goldschmidt's paper of 1974 \cite{goldschmidt74.45}.

\bq
{\normalsize
The author presents a new proof of the famous Brauer-Suzuki-Wall theorem: Let $G$ be a finite group and $T \in {\rm Syl}_2(G)$. Suppose $|T|=q>2$ and $T=C_G(t)$ for every $t\in T^\#$. Then either $T \trianglelefteq G$ or $G \simeq {\rm SL}(2,q)$. The highly elegant and very short proof not only avoids character theory completely but uses only a few elementary facts of group theory. In particular, the author demonstrates that the simple fact that two involutions always generate a dihedral group is a powerful tool in the above situation. Thus the structure of $M=N_G(T)$ is quickly determined and one is able to count the involutions in every coset of $M$ in $G$. At this point it is easy to determine the order of $G$. For the final identification of $G$ with ${\rm SL}(2,q)$ one shows that the permutation representation of $G$ on the cosets of $M$ is unique.
}
\eq

This was a slightly stronger result than Burnside's. Unfortunately, an obscure journal where it was published appeared to be unaccessible---but, obviously, Goldschmidt's paper was build around the same circle of ideas.\footnote{Otto Kegel wrote to me:
\bq
{\footnotesize From the left-overs of the three originally different papers that merged in the Brauer-Suzuki-Wall paper, Tim Wall and I wrote [in 19]61 a small paper
on groups with a partions \cite{kegel-wall61.255}.
}\\
\eq
This paper is interesting because it also avoids the use of theory of characters. It is a pity---I would love to know about it in 1982.}

But now I have Burnside's paper \cite{burnside00.269} in my hands, and I can see what was his thinking. The language was somewhat old-fashioned (he was saying ``operations'' instead of ``elements'' and ``self-conjugate'' subgroups instead of ``normal'' subgroups), but very clear---and remarkably in the spirit of Schmidt's book.

First of all, Burnside clearly uses Brauer's lemma on conjugacy of involutions:
\bq
{\normalsize
The operation $AB$ must
either be of order two or of odd order. If it were of odd order, $\mu$, the subgroup
generated by $A$ and $B$ would be a dihedral subgroup of order $2\mu$ ; and in this subgroup
$A$ and $B$ would be conjugate operations. [p. 271]
}
\eq
Then he observes that if involutions form more than one  conjugacy class (``conjugate set'', in his terminology), the subgroup generated by these classes commute pairwise, and it easily follows that a $2$-Sylow subgroup $H$ of $G$ is normal in $G$. Burnside analyses the structure of $G$ in that special case and then moves to
\bq
{\normalsize
the case in which the operations of $G$ of order two form a single conjugate set, while $G$ contains more than one subgroup of order $2^n$ [that is, $2$-Sylows---AB.].
}
\eq
He observes that $2$-Sylow subgroups are TI:
\bq
{\normalsize
no two subgroups of order $2^n$ have common operations other than identity.
}
\eq
On p. 272, Burnside introduces the normaliser $K = N_G(H)$ of $H$:
\bq
{\normalsize
\dots $K$ is the greatest subgroup of $G$ which contains a subgroup $H$, of order $2^n$,
self-conjugately \dots
}
\eq
and uses the previous analysis to determine the structure of $K$:
\bq
{\normalsize
\dots then $K$ must be a subgroup of the nature of those considered in the
preceding section, and its order must be $2^n\mu$, where $\mu$ is equal to or is a factor of $2^n - 1$.
}
\eq
In the next sentence, he formulates what is now called ``control of fusion'':
\bq
{\normalsize
Also no two operations of $H$ can be conjugate in $G$ unless they are conjugate
in $K$\footnote{\emph{Theory of Groups}, p. 98. [Burnside's footnote---AB.]}
}
\eq
and uses it to complete the analysis of the structure of $K$:
\bq
{\normalsize
The $2^n- 1$ operations of order two in $K$ therefore form a single conjugate set ;
and hence $\mu$ must be equal to $2^n -1$.
}
\eq
I skip further increasingly spooky details; by the end of p. 274 Burnside reaches the situation when $G$ is simple, has order $(2^n+1)\cdot 2^n\cdot(2^n-1)$, and acts faithfully and $3$-transitively on $2^n+1$ points.

Burnside's next statement---and the reader would perhaps agree with me---is something of a revelation:
\bq
{\normalsize
Hence
for a given value of $n$ the group, if it exists, is unique. [p.275]
}
\eq
I skip discussion of his proof of uniqueness; but what remains for Burnside to do after this has been done is to exhibit an example.
\bq
{\normalsize
That such groups exist for all values of $n$ is known\footnote{Moore : ``On a doubly-infinite series of simple groups,'' \emph{Chicago Congress Papers} (1893) ; Burnside : ``On a class of groups defined by congruences,'' \emph{Proc. L. M. S.} Vol. xxv. (1894). [Burnside's footnote---AB.]}. In fact the system of
congruences
\[
z' \equiv \frac{\alpha z + \beta}{\gamma z+ \delta}, \quad ({\rm mod}.\;\; 2),
\]
where $\alpha$, $\beta$, $\gamma$, $\delta$ are roots of the congruence
\[
\lambda^{2^n-1} \equiv 1, \quad ({\rm mod}.\;\; 2),
\]
such that
\[
\alpha\delta - \beta\gamma \not\equiv 0, \quad ({\rm mod}.\;\; 2),
\]
actually define such a group ; and the permutations of the $2^n + 1$ symbols

\[
\infty,\; 0,\; \lambda,\; \lambda^2, \dots, \lambda^{2^n-1},
\]
where $\lambda$ is a primitive root of
\[
\lambda^{2^n-1} \equiv 1, \quad ({\rm mod}.\;\; 2),
\]
which are effected by the above system of congruences, actually represent it as a triply-transitive group of degree $2^n+1$.
}
\eq

I hope you got the picture: William Burnside was working within the conceptual framework that much later became the classical paradigm of the classification of finite simple groups in terms of involutions.

Bill Kantor asked me on reading an earlier version of this manuscript:

\bq
{\normalsize
\dots what he [Burnside] did that influenced you so strongly?
}
\eq

What can I say? Burnside existed. More specifically, I became aware of him demonstrating that classification of objects similar and parallel to finite simple groups could be perhaps be carried out starting from the first principles. There was no need for me to read his proof, I could do it myself after spending years absorbing the group theoretic culture. But what I needed is to be told that this was possible and was worth trying.

\section{Euler, Hilbert, and others}
\label{sec:Euler}

But let me return my story back to the early 1980s and my first attempts to study groups of finite Morley rank. As I had already said, the key issue for me was understanding why  $G={\rm SO}_3(\mathbb{R})$ was not an $\omega$-stable group of finite Morley rank. The most natural way to achieve that was to interpret the field $\mathbb{R}$ in $G$ considered as an abstract group. If $G$ was $\omega$-stable and of finite Morley rank, then so would be the field $\mathbb{R}$, contradicting the celebrated theorem by Angus Macintyre: infinite $\omega$-stable fields are algebraically closed \cite{macintyre71.1}.

At that point I realised that I knew \emph{nothing} about ${\rm SO}_3(\mathbb{R})$, despite it being one of the most fundamental objects of mathematics\footnote{A naive question: why is the group ${\rm SO}_3(\mathbb{R})$ simple, while ${\rm SO}_3(\mathbb{Q})$ is not?  An even more fundamental fact: the cross product of vectors in $\mathbb{R}^3$ is the Lie algebra of ${\rm SO}_3(\mathbb{R})$---how many mathematics graduates are aware of that? Alas, fundamental concrete objects in mathematics are ignored in standard university mathematics courses.}. Luckily, I vaguely remembered that my university course of mechanics  once briefly mentioned \emph{the Euler angles}. I decided to try them and get a parametrisation of rotations by something which could be turn into a field. But, being within the configuration of the Brauer-Suzuki-Wall theorem as applied to ${\rm SO}_3(\mathbb{R})$, I had to start with involutions, that is, half-turns trough $\pi$ around axes in $\mathbb{R}^3$. Their crucial property is that two distinct half-turns commute if and only if their axes are perpendicular.

The ``conjugation by a square root'' trick came here very handy: if $i$ is an involution, $g \in G$ is an arbitrary element, take $j = i^g$ and use the formula
\[
j^{\sqrt{ji}} = i
\]
(keeping in mind that the square root is two-valued);
then
\[
g\cdot \sqrt{ji} = g\cdot \sqrt{i^g\cdot i} = g\cdot \sqrt{[g,i]} \in C_G(i).
\]
But $C_G(i)$ is the stabiliser in $G$ of the axis of $i$, and therefore its elements are associated with angles (I skip some minor details: how to handle the appearance of two roots, etc. -- this is fairly routine).

Taking three pairwise distinct and commuting involutions $i$, $j$, and $k$ is equivalent to choosing three perpendicular axes, and zeta-map associates with arbitrary rotation $g$ some three rotations about these axes and therefore three angles $\alpha, \beta, \gamma$ which parametrise ${\rm SO}_3(\mathbb{R})$.

But I had to overcome a fundamental difficulty of a  model-theoretic nature. Angles in the Euclidean space were measured by real numbers, but the reason for that was not of algebraic nature. In a purely algebraic set-up, angles were not expected to be numbers---I was aware of this from my school years thanks to the Russian translation of Dieudonn{\'e}'s book \cite{Dieudonne1963}\footnote{In 2005 I attended a fascinating conference on Euclid's \emph{Elements} in the Bodleian Library, Oxford. At one point historians of mathematics asked mathematicians, represented in the conference by Sir Christopher Zeeman and Robin Hartshorne, what aspect of \emph{Elements} they found to be the most important, and were quite surprised by the answer: Euclid never tried to measure angles by real numbers; angles and lengths were for him magnitudes of different nature.}. Therefore I had to work not with the angles of rotations, but with cosines of the angles, and the next obvious step was to use the fact that involutions commute if and only if their axes are perpendicular and try to develop some geometry of rectangular triangles. My first attempts to play with triangles quickly led me to the observation that the set of involutions in ${\rm SO}_3(\mathbb{R})$ attained a structure of a projective plane if  involution were taken for points and each set of involutions commuting with, but distinct from, the given involution $i$ called the line $l_i$; moreover, the maps $i \mapsto l_i$, $l_i \mapsto i$ formed a polarity of this plane. Moreover no kind of trigonometry appeared to be possible without first coordinatising this projective plane.

I started this story by telling about the little gem of mathematics exposition,  Hartshorne's book \cite{hartshorne67}. I knew from it Hilbert's famous coordinatisaion of Desargesian projective planes by  division rings. I had a projective plane, and I needed to make its Desargesian. I decided to invesigate where this plane came from,
and soon realised that the canonical representation of the rotation group ${\rm SO}_3(\mathbb{R})$ as the factor of the multiplicative group of the quaternions modulo the scalars
\[
{\rm SO}_3(\mathbb{R}) \simeq \mathbb{H}^*/\mathbb{R}^*
\]
naturally introduces on ${\rm SO}_3(\mathbb{R})$ the structure of a $3$-dimensional real projective space. Moreover, it was possible to prove within ${\rm SO}_3(\mathbb{R})$ seen as an abstract group that the plane made of involutions is indeed a plane in a $3$-dimensional projective space and is therefore Desarguisean, and, moreover, it was relatively easy to prove the commutativity of the resulting division ring. I got the field; the zeta map played in this construction the role of a temporary scaffolding and evaporated  without trace from the final version of the proof.

Alas, as soon as I mentioned to my colleagues my construction of the real projective plane within ${\rm SO}_3(\mathbb{R})$  given as an abstract group, I was immediately told that this construction had already been published by Bachmann in  his book \cite{bachmann1959} of 1959. Also, I found a beautiful paper by Boris Weisfeiler \cite{weisfeiler79.522} which explained the reasons for projective planes being interpretable in $k$-forms of ${\rm PSL}_2$.

But I think now that my further study of groups of finite Morley rank would be impossible without this amusing re-discovery of the geometry of involutions in  ${\rm SO}_3(\mathbb{R})$.

But at that time, my work in model-theoretic algebra had been interrupted.
Very soon after completion of my work on ${\rm SO}_3(\mathbb{R})$ I had to abandon my study groups of finite Morley rank for almost five years because this particular theme was distinctively unfashionable in Novosibirsk. Part of the story is told in \emph{Introduction} to \cite{borovik_ABC_2008}. I continue my earlier quote from \cite{borovik_ABC_2008}:

\bq
{\normalsize
{A year later Simon Thomas sent to me the manuscripts of
his work on locally finite groups of finite Morley rank. Besides
many interesting results and observations his manuscripts contained also an exposition of Boris Zilber's fundamental results on $\aleph_1$-categorical structures which were made known to many western model theorists in Wilfrid Hodges' translation of Zilber's paper \cite{zilber77viniti} but which, because of the regrettably restricted form of publication of the Russian original, remained unknown to me.}
}
\eq

Indeed I soon discovered that there were no way for me to publish my results in any form more serious than some obscure preprints \cite{borovik84.19,borovik84.21}\footnote{But even these preprints would not appear without an expression of interest and support from Otto Kegel who visited Novosibirsk in May 1983 and heard my talk on involutions at the impromptu ``Kegelfest'' in a coffee room of the Institute of Mathematics. (On receiving the news of his visit,  algebraists from all over Siberia rushed to Novosibirsk to give a talk on their latest research in front of Otto Kegel, who, slightly stunned by the unexpected burden but attentive and supportive, was patiently listening to everyone for three or four days in a row. But this Bulgakovesque episode deserves to be told in  a separate tale.)}. In these sad years not only the theme was unfashionable---the name of Zilber was unfashionable in Novosibirsk, too. The situation changed only when my work started to attract attention of model theorists in the West (and Otto Kegel's help was crucially important for establishing communication) and I started to develop collaboration with colleagues in the international mathematical community \cite{borovik90.478}. In particular, I used some ideas from my proof of Burnside's Theorem in my work with Ali Nesin on CIT groups of finite Morley rank \cite{borovik94.258,borovik94.273}.

My work on involutions in ${\rm SO}_3(\mathbb{R})$ triggered my interest to algebraic aspect of the theory symmetric spaces; I had even taught an advanced level university course based on books by Ottmar Loos \cite{loos69-I,loos69-II}. Not surprisingly, the name of \'Elie Cartan featured in the theory of symmetric spaces prominently. It will appear in my story again, see Section~\ref{sec:Cartan}.

\section{Black box groups}

But this is still not the end of the story: the  ``conjugation by a square root'' and ``double conjugation'' tricks had happened to be very useful in the black box group theory, too.

Many years later, in July 1996, Laci Babai told me about the difficulty of distinguishing in the black box setting between the symplectic group $C_n(q) = {\rm PSp}_{2n}(q)$ and
the orthogonal group $B_n(q) = \Omega_{2n+1}(q)$ where $q>3$ is odd  and $n$ and $q$ are given. Indeed, since these groups have the same Weyl group,  their statistics of orders of random elements are undistinguishable in polynomial time. I was able to tell Laci that these groups had completely different structures of centralisers of involutions and that centralisers of involutions could be computed in the black box setup. Our conversation was at a group theory conference in Pisa. On my return to Manchester I offered this problem to my PhD student Christine Altseimer who already worked with involutions in groups of finite Morley rank \cite{altseimer99.1879}. This eventually resulted in our paper \cite{altseimer01.1} which contained an answer to Babai's question.

The ``double conjugation'' trick is an established part of group-theoretical folklore; Bill Kantor, in particular, reminded me that it
\bq
{\normalsize
has been discovered and used elsewhere, certainly in black box settings in papers of mine.
}
\eq
In my own work, the current state of the ``double conjugation'' trick is the following simple observation which plays the crucial role in proofs of Theorems 3.1 and 5.2 of the paper \cite{suko12F} by \c{S}\"{u}kr\"{u} Yal\c{c}\i nkaya and myself:

\bq
{\normalsize
\textbf{The double conjugation trick.} Let $X$ be a black box group and $U<X$ an elementary abelian TI $2$-subgroup given as a black box subgroup of $X$.
Assume in addition that the action of $N=N_X(U)$  on $U^\#$ by conjugation is transitive.
Fix an involution $v \in U^\#$. Then we have a partial map
\[
\nu_v: U^\#  \times X \longrightarrow N_X(U)
\]
defined on all pairs $(u,x) \in U^\#  \times X$ such that  the element $u^x\cdot v$ has odd order:
\[
\nu_v: (u,x) \mapsto x\sqrt{u^x \cdot v}
\]
Moreover, if $n \in N$ then
\[
\nu_v(u^{-n}, nx) = nx \cdot\sqrt{(u^{-n})^{nx} \cdot v} =nx\sqrt{u^x \cdot v} = n\nu(u,x),
\]
which means that the probability distribution of elements $\nu_v(u,x)$ in $N$ is invariant under the left action by $N$. One may say that
\bq
{\normalsize
\emph{A recombination of black boxes for $U$ and $X$ produces a black box for $N=N_X(U)$.}
}
\eq
}
\eq

Notice that if $|U|=2$, the procedure described turns out to be Bray's algorithm as described in Section~\ref{sec:centralisers} and
\[
\zeta(x) = \nu_u(u,x) = x\sqrt{u^x \cdot u}
\]
is the corresponding partial  map
\[
\zeta: X \longrightarrow C_X(u).
\]
In (a pretty hypothetical) case of black box groups of even characteristic a canonical example of TI subgroups is given by long root subgroups and the double conjugation trick become the following theorem.

\bq
 {\normalsize
\textbf{Theorem 5.2 \cite{suko12F}} Let\/ $X$ be a black box group encrypting an untwisted Chevalley group ${\rm G}(2^n)$ {\rm (}with $n$ known{\rm )} and\/ $U < X$ an unipotent long root subgroup given as a black box subgroup of\/ $X$. Then there is a polynomial time, in $n$ and  Lie rank of ${\rm G}(2^n)$, Monte-Carlo algorithm which constructs a black box for $N_X(U)$ and a black box $\mathbb{U}$ for the field  $\mathbb{F}_{2^n}$ interpreted in the action of $N_X(U)$ on $U$, with $U$ becoming the additive group of the field $\mathbb{U}$.
 }
\eq

This result suggests that if by any chance we are given an involution in a black box group  $X$ encrypting a  Chevalley groups ${\rm G}(2^n)$ of characteristic $2$ (admittedly, a very big ``if''), then its structure recovery  is likely to share some of the conceptual framework of Franz Timmesfeld's classification of groups generated by root type subgroups \cite{timmesfeld90.167,timmesfeld91.575,timmesfeld92.183,timmesfeld2001}. If so, then this will be strikingly similar to the use of Aschbacher's classical involutions \cite{aschbacher77.353,aschbacher80.411} and root ${\rm SL}_2$-subgroups in the structural theory of classical black box groups in odd characteristic \cite{suko03,suko12E,suko12A,suko12B,suko02,yalcinkaya07.171} as developed by \c{S}\"{u}kr\"{u} Yal\c{c}\i nkaya and myself.

I have to postpone  to another occasion a detailed discussion of parallels between the recursive algorithms for black box groups recognition and the inductive proof of the classification of finite simple groups and especially  black box analogues of Aschbachers' Classical Involution Theorem as first discovered by Christine Altseimer and then developed by \c{S}\"{u}kr\"{u} Yal\c{c}\i nkaya into a flexible an powerful method not only of identification, but detailed and intricate analysis of black box groups One of the most interesting aspects of it is the role and nature of the so-called ``signaliser functors'' and ``$L$-balance'' so prominent in the classification of finite simple groups and appearing in two different but complimentary aspects in the theory of groups of finite Morley rank (Ay\c{s}e Berkman and Jeff Burdges) and in the black box group theory (\c{S}\"{u}kr\"{u} Yal\c{c}\i nkaya).

\section{\'Elie Cartan and others}
\label{sec:Cartan}

In the specific case when $t \in X$ is a strongly isolated involution, that is, when $t$ is not conjugate  in $X$ to any other involution in a $2$-Sylow subgroup $S$ containing $t$, elements
\[
t^xt = x^{-1}\cdot x^t
\]
have odd order for all $x\in X$ and the map
\[
\zeta(x) = x\sqrt{x^{-1}\cdot x^t}
\]
is defined for all $x\in X$ and is $C_X(t)$-equivariant: for all $c\in C_X(t)$ we have
\bea
\zeta (cx) &=& cx\cdot \sqrt{(cx)^{-1}\cdot (cx)^\tau}\\
    &=& cx\cdot \sqrt{x^{-1}c^{-1}c^\tau x^\tau}\\
    &=& cx\cdot\sqrt{x^{-1}c^{-1}c x^\tau}\\
    &=& cx\cdot\sqrt{x^{-1} x^\tau}\\
    &=& c\zeta(x).
\eea
A finite group $X$ with a strongly isolated involution $t$ cannot be simple, this is Glauberman's celebrated $Z^*$-Theorem \cite{glauberman66.403}:
\[
t \in Z^*(X) = Z(X \bmod O(X)).
\]
An analogue of Glauberman's $Z^*$-Theorem for groups $X$ of finite Morley rank remains a widely open problem, but we still have a definable $C_X(t)$-equivariant map
\[
\zeta: X \longrightarrow C_X(t).
\]
In particular, if $X$ is connected, then the centraliser $C_X(t)$ of a strongly isolated involution $t$ is also connected. This application of black box methods to groups of finite Morley rank has a remarkable range of profound consequences \cite{borovik07.1}.

When I first mentioned this observation to my model theorists friends, I was immediately told by  Adrien Deloro that the same trick can be used  to prove connectedness of the real orthogonal group ${\rm SO}_n(\mathbb{R})$: it is the centraliser in ${\rm SL}_n (\mathbb {R})$ of the inverse-transposed automorphism $\tau: x \mapsto (x^{-1})^T$, where $T$ is the matrix transposition; the associated map
\bea
\zeta: {\rm SL}_n (\mathbb {R}) & \longrightarrow &  {\rm SO}_n(\mathbb{R})\\
 x & \mapsto & x\sqrt{x^{-1}\cdot x^\tau}
\eea
is defined for all $x \in {\rm SL}_n (\mathbb {R})$ because the matrix
\[
x^{-1}\cdot x^\tau = x^{-1}\cdot (x^{-1})^T
\]
is a positively defined symmetric matrix, hence is orthogonally diagonalisable with real positive eigenvalues and thus has in ${\rm SL}_n (\mathbb {R})$ a (uniquely defined) square root with positive eigenvalues.
Moreover, $\zeta$ is continuous and ${\rm SO}_n(\mathbb{R})$-equivariant,
\[
\zeta(sx) = s\zeta(x) \mbox{ for all } s\in {\rm SO}_n(\mathbb{R}),
\]
therefore the connectedness of ${\rm SL}_n (\mathbb {R})$ implies the connectedness of ${\rm SO}_n(\mathbb{R})$.

If we denote
\[
z = x \sqrt{x^{-1}\cdot (x^{-1})^T}
\]
we get the classic \emph{polar decomposition}
\[
x = z \cdot \sqrt{x^T\cdot x}
\]
of $x$ as a product of an orthogonal matrix $z$ and the positively defined symmetric matrix $\sqrt{x^T\cdot x}$. Similarly, for matrices  $x\in {\rm GL}_n(\mathbb{C})$ we have the polar decomposition
\[
x = z \cdot \sqrt{\bar{x}^T\cdot x}
\]
of $x$ into the product of an unitary matrix $z$ and the positively determined hermitian matrix $\sqrt{\bar{x}^T\cdot x}$.

What Adrien Deloro told me next was even more striking. In his recent email he has refreshed for me the details of the story:

\bq
{\normalsize
At classes pr\'{e}pa\footnote{\cite{lemme12.5} contains an informative discussion of \emph{classes pr\'{e}pa}.}  level the polar decomposition is sometimes called \emph{Cartan's
decomposition}. There are of course reasons to that (the proximity both of
Cartan's decomposition of a Lie \emph{algebra} and of the $KAK$ decomposition---called Cartan's decomposition in more serious sources).
As a result the terminology is rather systematically abused in France.

But I can't find a single \emph{textbook} calling the polar decomposition Cartan's.
Only \emph{pr\'{e}pa} level \emph{lecture notes} do so, and you are right saying the confusion is oral tradition.

I do not know whether \'{E}lie Cartan used the polar decomposition to prove
connectedness; what I know is that this method has become ``the standard
example'' at an \emph{agr\'{e}gation} lecture.
}
\eq

Indeed, the authoritative monograph by  Horn and Johnson \cite{horn-johnson1902} attributes the discovery of polar decomposition of square complex matrices to Autonne and his  paper \cite{autonne1902} of 1902; the authors do not mention Cartan at all---perhaps because their book is written from the viewpoint of functional analysis, not algebra. But, as Adrien Deloro rightly emphasises, the polar decomposition is a special case of the \emph{Cartan decomposition} which works in all semi-simple complex Lie algebras and Lie groups; so perhaps Cartan was involved with more elementary cases of his decomposition, as well.

Still, at this historical background \c{S}\"{u}kr\"{u} Yal\c{c}\i nkaya and I are pleased that, given a black box group $X$ encrypting ${\rm SL}_n(q)$ for odd $q$,  we have constructed the inverse-transpose automorphism $\tau$ and therefore got an algorithm for polar decomposition:
 \[
x = x^\tau \cdot (x^\tau)^{-1}x,
 \]
 and then used the zeta map (which in this context could be rightly called Cartan's map) to construct the subgroup ${\rm SO}_n(q)$ in $X$ as the centraliser $C_X(\tau)$ of the involution $\tau$ \cite{suko12B}. Notice that $\tau$ was not part of the original setup of the black box $X$!

And the final word on the role of tradition belongs to Adrien Deloro:\label{deloro-quote}
 \bq
 {\normalsize
 As you may know, I never lost hope to prove something on the subject---in the finite Morley rank context, in odd type. I always wanted to dedicate it to my \emph{pr\'{e}pa} teacher. The bitter irony is, three years ago I had in my personal notes started denoting your zeta map capital $\acute{\mathcal{E}}$,  for \'{E}lie (which in \texttt{$\backslash$mathcal} also looks like zeta). Of course \'{E}ric's insisting that we ought to mention the zeta map in the locally solvable case (our paper \cite{Deloro-Jaligot10.23} at Crelle's journal), and his death\footnote{Our colleague and friend \'{E}ric Jaligot died on 12 July 2013 at the peak of his mathematical creativity.}  made such a tribute impossible for another couple of years \dots
 }
 \eq

\section{Conclusion}

This paper documents a convoluted history of a few simple but powerful mathematical ideas and the way they were inherited, or ignored, or re-discovered by new generations of mathematicians.

I was privileged to work with Israel Moiseevich Gelfand. He made a clear distinction between the two modes of work in mathematics expressed by Russian words  \emph{pr\textbf{I}dumyvanie} and \emph{pr\textbf{O}dumyvanie}, very similar and almost homophonic. The former means `inventing', the latter `properly thinking through' and was used by Gelfand with the meaning
 \bq
 `\emph{systematically thinking through back to the origins, fundamentals, first principles of particular mathematical concepts or problems.}'
 \eq
Gelfand valued \emph{produmyvanie} more than \emph{pridumyvanie}, he believed that \emph{produmyvanie} yielded deeper results. The stories told here, I hope, support this view.

There was also another reason for me to spend  time on documenting this story. On a number of occasions I was asked by my philosopher colleagues about possibility of a ``Lakatosian'' mathematical education \cite{larvor10.71} built around students' inquiry in a free flowing dialectic development of mathematical concepts, results and proofs.

I belong to a peculiar group of people whose education was a crude approximation to this ideal. As many Russian mathematicians of my generation, I was a product of the ``olympiad'' system, a loose network of mathematics competitions, outreach activities, correspondence courses for school children, and specialist boarding schools run by elite universities.

The system was disorganised, chaotic, and unaccountable. But it worked, and this paper provides a glimpse into its inner working and shows that it had a distinctive Lakatosian flavour; it appeared because mathematics education of school children  was handled by professional research mathematicians. The interested  reader may wish to see \cite{borovik12.23} for a brief sketch of life and study at my \emph{alma mater},  the FMSh, the Preparatory Boarding School of the Novosibirsk University. I hope that \cite{borovik12.23} closes the question of feasibility of a ``Lakatosian'' model of mathematics education in the real world: as all really good things, it is economically unsustainable.

\section*{Acknowledgements}

I will be forever grateful to Eric Jaligot (1972--2013) for his enthusiastic interest to the zeta map; he had even called it ``the magic map''.

I am grateful to Bill Kantor for his very stimulating emails, to Adrien Deloro and \c{S}\"{u}kr\"{u} Yal\c{c}\i nkaya for sharing with me the joys and sorrows of mathematical research, and to Gregory Cherlin and Otto Kegel for their support and help over many years---and for useful comments on this paper.

\section*{Disclaimer}

\bq
\emph{The author writes in his personal capacity and  the views expressed do not necessarily represent position of his employer or any other organisation or institution.}
\eq

\bibliographystyle{amsplain}
\bibliography{lacatosian}

\def\cprime{$'$}
\providecommand{\bysame}{\leavevmode\hbox to3em{\hrulefill}\thinspace}
\providecommand{\MR}{\relax\ifhmode\unskip\space\fi MR }
\providecommand{\MRhref}[2]{%
  \href{http://www.ams.org/mathscinet-getitem?mr=#1}{#2}
}
\providecommand{\href}[2]{#2}
\begin{thebibliography}{10}

\bibitem{borovik_ABC_2008}
T.~Alt{\i}nel, A.~Borovik, and G.~Cherlin, \emph{Simple groups of finite
  {M}orley rank}, Mathematical {M}onographs, Amer. Math. Soc., 2008.

\bibitem{altseimer99.1879}
C.~Altseimer, \emph{A characterisation of {${\rm PSp}(4,K)$}}, Comm. Algebra
  \textbf{27} (1999), no.~4, 1879--1888. \MR{1679699 (2000a:20070)}

\bibitem{altseimer01.1}
C.~Altseimer and A.~V. Borovik, \emph{Probabilistic recognition of orthogonal
  and symplectic groups}, Groups and computation, {III} ({C}olumbus, {OH},
  1999), Ohio State Univ. Math. Res. Inst. Publ., vol.~8, de Gruyter, Berlin,
  2001, pp.~1--20; Corrections: math.GR/0110234. \MR{1829468 (2002e:20093)}

\bibitem{aschbacher77.353}
M.~Aschbacher, \emph{A characterization of {C}hevalley groups over fields of
  odd order. {I, II}}, Ann. of Math. (2) \textbf{106} (1977), no.~3, 353--468.

\bibitem{aschbacher80.411}
\bysame, \emph{Correction to: ``{A} characterization of {C}hevalley groups over
  fields of odd order. {I}, {II}'' [{A}nn. of {M}ath. (2) {\bf 106} (1977),
  353--468]}, Ann. of Math. (2) \textbf{111} (1980), no.~2, 411--414.

\bibitem{autonne1902}
L.~Autonne, \emph{Sur le groupes lin\'{e}ares, r\'{e}els et orthogonaux}, Bull.
  Soc. Math. France \textbf{30} (1902), 121--134.

\bibitem{bachmann1959}
F.~Bachmann, \emph{Aufbau der {G}eometrie aus dem {S}piegelungsbegriff},
  Springer-Verlag, Berlin-G\"ottingen-Heidelberg, 1959.

\bibitem{borovik07.1}
A.~Borovik, J.~Burdges, and G.~Cherlin, \emph{Involutions in groups of finite
  {M}orley rank of degenerate type}, Selecta Math. (N.S.) \textbf{13} (2007),
  no.~1, 1--22. \MR{2330585 (2008f:03061)}

\bibitem{borovik84.19}
A.~V. Borovik, \emph{Involyutsii v gruppakh s razmernostyu}, Preprint, vol.~84,
  Akad. Nauk SSSR Sibirsk. Otdel. Vychisl. Tsentr, Novosibirsk, 1984.
  \MR{810216 (87j:20071b)}

\bibitem{borovik84.21}
\bysame, \emph{Teoriya konechnykh grupp i neschetno kategorichnye gruppy},
  Preprint, vol.~84, Akad. Nauk SSSR Sibirsk. Otdel. Vychisl. Tsentr,
  Novosibirsk, 1984. \MR{818985 (87j:20071a)}

\bibitem{borovik02.7}
\bysame, \emph{Centralisers of involutions in black box groups}, Computational
  and statistical group theory (Las Vegas, NV/Hoboken, NJ, 2001), Contemp.
  Math., vol. 298, Amer. Math. Soc., Providence, RI, 2002, pp.~7--20.

\bibitem{borovik12.23}
\bysame, \emph{Free maths schools: some international parallels}, The De Morgan
  J. \textbf{2} (2012), no.~2, 23--35,
  http://education.lms.ac.uk/wp-content/uploads/2012/02/FMSh.pdf.

\bibitem{borovik94.258}
A.~V. Borovik, M.~J. DeBonis, and A.~Nesin, \emph{C{IT} groups of finite
  {M}orley rank. {I}}, J. Algebra \textbf{165} (1994), no.~2, 258--272.
  \MR{1273275 (95h:03082)}

\bibitem{borovik94.273}
A.~V. Borovik and A.~Nesin, \emph{C{IT} groups of finite {M}orley rank. {II}},
  J. Algebra \textbf{165} (1994), no.~2, 273--294. \MR{1273276 (95h:03083)}

\bibitem{borovik90.478}
A.~V. Borovik and B.~P. Poizat, \emph{Tores et {$p$}-groupes}, J. Symbolic
  Logic \textbf{55} (1990), no.~2, 478--491. \MR{1056365 (91j:03045)}

\bibitem{suko03}
A.~V. Borovik and {\c{S}}.~Yal\c{c}{\i}nkaya, \emph{Construction of
  {C}urtis-{P}han-{T}its system for black box classical groups}, Available at
  arXiv:1008.2823v1 [math.GR].

\bibitem{suko12A}
\bysame, \emph{Steinberg presentations of black box classical groups in small
  characteristics}, Available at arXiv:1302.3059v1 [math.GR].

\bibitem{suko12F}
\bysame, \emph{Fifty shades of black}, Available at arXiv:1308.2487 [math.GR].

\bibitem{suko12E}
\bysame, \emph{Classical black box groups in small odd characteristics}, in
  preparation.

\bibitem{suko12B}
\bysame, \emph{Subgroup structure and automorphisms of black box classical
  groups}, in preparation.

\bibitem{brauer58.718}
R.~Brauer, M.~Suzuki, and G.~E. Wall, \emph{A characterization of the
  one-dimensional unimodular projective groups over finite fields}, Illinois J.
  Math. \textbf{2} (1958), no.~4B, 718--745.

\bibitem{bray00.241}
J.~N. Bray, \emph{An improved method for generating the centralizer of an
  involution}, Arch. Math. (Basel) \textbf{74} (2000), no.~4, 241--245.

\bibitem{burnside00.269}
W.~Burnside, \emph{On a class of groups of finite order}, Transactions of the
  Cambridge Philos. Soc. \textbf{18} (1900), 269--276.

\bibitem{burnside1911}
\bysame, \emph{Theory of groups of finite order}, Cambridge Univ. Press, London
  and New York, 1911.

\bibitem{cherlin1979.1}
G.~Cherlin, \emph{Groups of small {M}orley rank}, Ann. Math. Logic \textbf{17}
  (1979), no.~1, 1--28.

\bibitem{Deloro-Jaligot10.23}
A.~Deloro and Eric Jaligot, \emph{{Small groups of finite Morley rank with
  involutions}}, Journal Fur Die Reine Und Angewandte Mathematik \textbf{2010}
  (2010), 23--45.

\bibitem{Dieudonne1963}
J.~Dieudonn{\'e}, \emph{Alg\`ebre lin\'eaire et g\'eom\'etrie \'el\'ementaire},
  Enseignement des Sciences, VIII, Hermann, Paris, 1964. \MR{0171788 (30
  \#2015)}

\bibitem{feit79.1}
W.~Feit, \emph{Richard {D}. {B}rauer}, Bull. Amer. Math. Soc. (New Ser.)
  \textbf{1} (1979), no.~1, 1--20.

\bibitem{glauberman66.403}
G.~Glauberman, \emph{Central elements in core-free groups}, J. Algebra
  \textbf{4:} (1966), 403--420.

\bibitem{goldschmidt74.45}
D.~M. Goldschmidt, \emph{Elements of order two in finite groups}, Delta
  (Waukesha) \textbf{4} (1974/75), 45--58. \MR{0352256 (50 \#4743)}

\bibitem{gorenstein1968}
D.~Gorenstein, \emph{Finite groups}, second ed., Chelsea Publishing Co., New
  York, 1980. \MR{81b:20002}

\bibitem{hall-higman56.1}
P.~Hall and G.~Higman, \emph{On the {$p$}-length of {$p$}-soluble groups and
  reduction theorems for {B}urnside's problem}, Proc. London Math. Soc. (3)
  \textbf{6} (1956), 1--42. \MR{0072872 (17,344b)}

\bibitem{hartshorne67}
R.~Hartshorne, \emph{Foundations of projective geometry}, Lecture Notes,
  Harvard University, vol. 1966/67, W. A. Benjamin, Inc., New York, 1967.
  \MR{0222751 (36 \#5801)}

\bibitem{horn-johnson1902}
R.~A. Horn and C.~R. Johnson, \emph{Topics in matrix analysis}, Cambridge
  University Press, Cambridge, 1991. \MR{1091716 (92e:15003)}

\bibitem{kegel-wall61.255}
O.~H. Kegel and G.~E. Wall, \emph{Zur {S}truktur endlicher {G}ruppen mit
  nicht-trivialer {P}artition}, Arch. Math. (Basel) \textbf{12} (1961),
  255--261. \MR{0136648 (25 \#114)}

\bibitem{larvor10.71}
B.~Larvor, \emph{Authoritarian versus authoritative teaching: {P}olya and
  {L}akatos}, Explanation and Proof in Mathematics (G.~Hanna, H.~N. Jahnke, and
  H.~Pulte, eds.), Springer US, New York, 2010, pp.~71--83.

\bibitem{lemme12.5}
M.~Lemme, \emph{Utter elitism: French mathematics and the system of classes
  pr\'{e}pas}, The De Morgan J. \textbf{2} (2012), no.~2, 5--22,
  \url{http://education.lms.ac.uk/wp-content/uploads/2012/02/Louis-le-Grand1.pdf}.

\bibitem{loos69-II}
O.~Loos, \emph{Symmetric spaces. {I}: {G}eneral theory}, W. A. Benjamin, Inc.,
  New York-Amsterdam, 1969. \MR{0239005 (39 \#365a)}

\bibitem{loos69-I}
\bysame, \emph{Symmetric spaces. {II}: {C}ompact spaces and classification}, W.
  A. Benjamin, Inc., New York-Amsterdam, 1969. \MR{0239006 (39 \#365b)}

\bibitem{macintyre71.1}
A.~Macintyre, \emph{On $\omega_1$-categorical theories of fields}, Fundamenta
  Mathematicae \textbf{71} (1971), 1--25.

\bibitem{matveeva93}
L.~V. Matveeva, \emph{Otto {Y}ul'evich {S}hmidt. 1891--1956},
  {Nauchno-Biograficheskaya Seriya}. [Scientific-Biographic Series], ``Nauka'',
  Moscow, 1993. \MR{1270211 (95k:01026)}

\bibitem{schmidt66}
O.~U. Schmidt, \emph{Abstract theory of groups}, Translated from the Russian by
  Fred Holling and J. B. Roberts. Translation edited by J. B. Roberts, W. H.
  Freeman and Co., San Francisco, Calif., 1966. \MR{0209340 (35 \#238)}

\bibitem{timmesfeld90.167}
F.~G. Timmesfeld, \emph{Groups generated by {$k$}-transvections}, Invent. Math.
  \textbf{100} (1990), no.~1, 167--206.

\bibitem{timmesfeld91.575}
\bysame, \emph{Groups generated by {$k$}-root subgroups}, Invent. Math.
  \textbf{106} (1991), no.~3, 575--666.

\bibitem{timmesfeld92.183}
\bysame, \emph{Groups generated by {$k$}-root subgroups---a survey}, Groups,
  combinatorics \& geometry (Durham, 1990), London Math. Soc. Lecture Note
  Ser., vol. 165, Cambridge Univ. Press, Cambridge, 1992, pp.~183--204.

\bibitem{timmesfeld2001}
\bysame, \emph{Abstract root subgroups and simple groups of {L}ie type},
  Monographs in Mathematics, vol.~95, Birkh\"{a}user Verlag, Basel, 2001.

\bibitem{weisfeiler79.522}
B.~Weisfeiler.

\bibitem{suko02}
{\c{S}}.~Yal\c{c}{\i}nkaya, \emph{Construction of long root
  {S}{L}$_2(q)$-subgroups in black-box groups}, Available at arXiv,
  math.GR/1001.3184v1.

\bibitem{yalcinkaya07.171}
\c{S}. Yal{\c{c}}{\i}nkaya, \emph{Black box groups}, Turkish J. Math.
  \textbf{31} (2007), no.~suppl., 171--210. \MR{2369830 (2009a:20081)}

\bibitem{zilber77viniti}
B.~I. Zilber, \emph{The structure of models of categorical theories and the
  problem of finite axiomatizability}, Manuscript deposited at VINITI on July
  12, 1977, Deposition No. 2800--77, 1977.

\end{thebibliography}

\small
\textsc{Email} \texttt{alexandre.borovik \(>>>\)at\(<<<\) gmail.com}\\

\end{document}